\documentclass[11pt]{article}

\usepackage{amscd}  
\usepackage{amsfonts}  
\usepackage{amssymb}  
\usepackage{latexsym}  
 
\newcommand{\sect}[1]{\section{#1}\setcounter{equation}{0}}

\newtheorem{definition}{Definition}[section]

\newtheorem{example}[definition]{Example}

\newtheorem{proposition}[definition]{Proposition}
\newtheorem{theorem}[definition]{Theorem}

\newtheorem{lemma}[definition]{Lemma}

\def\Gal{\mathop{\rm Gal}\nolimits}

\def\Alg{\mathsf{Alg}}
\def\kAlg{\mbox{$k$-$\mathsf{Alg}$}} 
\def\alg{\mbox{-}\mathsf{alg}}
 
\def\MonCat{\mathsf{MonCat}}
\def\Bia{\mathsf{Bia}}

\def\M{\mathcal{M}} 
\def\G{\mathcal{G}}

\def\o{\otimes}
\def\x{\times}
\def\rcross{\rtimes}
\def\bra{\langle} 
\def\ket{\rangle} 
\def\under{\mbox{\rm\_}\,} 
\newcommand{\rarr}[1]{\stackrel{#1}{\longrightarrow}} 
\newcommand{\larr}[1]{\stackrel{#1}{\longleftarrow}} 
\def\iso{\rarr{\sim}}
\def\cop{\Delta}

\def\eps{\varepsilon} 
 
\def\du1{\hat 1}

\newcommand{\End}{\mbox{\rm End}}
\newcommand{\Aut}{\mbox{\rm Aut}}
\newcommand{\Hom}{\mbox{\rm Hom}}
\newcommand{\Center}{\mbox{\rm Center}}
\def\id{\mbox{\rm id}}
\def\Ad{\mbox{\rm Ad}}
\def\0{_{(0)}} 
\def\1{_{(1)}} 
\def\2{_{(2)}} 
\def\3{_{(3)}} 
\def\PL{\pi^{\scriptscriptstyle L}} 
\def\L{^{\scriptscriptstyle L}} 
\def\R{^{\scriptscriptstyle R}} 
\def\PR{\pi^{\scriptscriptstyle R}}

\def\lact{\triangleright} 
 
\def\op{^{\rm op}}
\newcommand{\Sub}{\mbox{\rm Sub}}
\newcommand{\Fix}{\mbox{\rm Fix}}

\begin{document}

\title{Galois actions by finite quantum groupoids}

\author{\sc Korn\'el Szlach\'anyi\\
Research Institute for Particle and Nuclear Physics\\
H-1525 Budapest, P.O.Box 49\\
email:\,\tt szlach@rmki.kfki.hu
}
\date{}

\maketitle

\footnotetext{To appear in the
Proceedings of the "69$^{\mbox{\tiny \`eme}}$ rencontre entre physiciens
th\'eoriciens et math\'ematiciens", Strasbourg 2002}
\footnotetext{Partially supported by Hungarian Scientific Research Fund,
OTKA T-034512}

\vskip 1.6truecm

\begin{abstract}
Proposing a certain category of bialgebroid maps we show that the balanced
depth 2 extensions appear as they were the finitary Galois extensions in the
context of quantum groupoid actions, i.e., actions by finite bialgebroids, weak
bialgebras or weak Hopf algebras. We comment on deformation of weak
bialgebras, on half grouplike elements, on uniqueness of weak Hopf algebra
reconstructions and discuss the example of separable field extensions.
\end{abstract}

\vskip 1.2truecm

For extensions of rings, algebras or $C^*$-algebras the notion of depth 2,
introduced originally for von Neumann factors by A. Ocneanu, has many features
that makes it the analogue of Galois extension of fields. The extension
$N\subset M$ of $k$-algebras is called of depth 2 if the canonical $N$-$M$
bimodule $X=\,_NM_M$ and $M$-$N$ bimodule $\bar X=\,_MM_N$ satisfy: $X\o\bar
X\o X$ is a direct summand in a finite direct sum of copies of $X$ and $\bar
X\o X\o\bar X$ is a direct summand in a finite direct sum of $\bar X$'s. The
right module $M_N$ is called balanced if $\End\,_\mathcal{E}M\cong N$ where
$\mathcal{E}=\End\,M_N$. For any balanced depth 2 extension the endomorphism
ring $A=\End\,_NM_N$ carries a bialgebroid structure which is finite
projective over the centralizer, or relative commutant, $R=C_M(N)$ both as a
left and as a right module. Moreover the canonical action of $A$ on $M$ makes
$M$ to be a left $A$-module algebroid with invariant subalgebra equal to $N$
\cite{KSz}. This generalizes the result of Nikshych and Vainerman in subfactor
theory \cite{NV1} because finite index depth 2 extensions of von Neumann or
$C^*$-algebras are always balanced depth 2 extensions. Moreover they are
Frobenius extensions which causes the appearence of antipodes, hence leading
to weak $C^*$-Hopf algebra structure \cite{BSz} on $A$.

Our main purpose in this paper will be to study the uniqueness problem of the
bialgebroid $A$. Using a kind of category for bialgebroids we show that $A$
satisfies a universal property analogous to the one of the Galois group of a
field extension. While the natural action of $A$ on $M$ generalizes
Hopf-Galois extensions of Kreimer and Takeuchi \cite{Kr-Ta}, Hopf-Galois
extensions do not have the universal property w.r.t the category of Hopf
algebras as the example of Greither and Pareigis \cite{Gr-Pa} has shown. 
Our proposal of a Galois bialgebroid $\Gal(M/N)$ of an algebra extension in
Section \ref{sec: Gal} is close in spirit to Pareigis' Quantum Automorphism
Group \cite{Pareigis} but technically not a mature one. 

In case of depth 2 Frobenius extensions we present an optimistic
interpretation of the non-uniqueness of its associated weak Hopf algebra. The
natural object to which a unique (measured) weak Hopf algebra can be
associated is a Frobenius system \cite{K}, i.e., a Frobenius extension
$N\subset M$ together with a Frobenius homomorphism
$\phi\colon\,_NM_N\to\,_NN_N$. 

Let me recall \cite{BNSz} for the definition of weak bialgebra (WBA) which is
one of the main theme in this paper. Let $K$ be a field. A finite dimensional
$K$-space $A$ together with a $K$-algebra structure $\bra
A,m,u\ket$ and a $K$-coalgebra structure $\bra A,\cop,\eps\ket$ is called a
\textbf{weak bialgebra} if

\begin{enumerate}
\item $\cop$ is multiplicative / $m$ is comultiplicative, i.e., as maps $A\o
A\to A\o A$, 
\[
\cop\circ m=(m\o m)\circ(\id \o\Sigma\o\id )\circ(\cop\o\cop)
\]
where $\Sigma\colon A\o A\to A\o A$ denotes the flip map $x\o y\mapsto y\o x$,
\item $\eps$ is weakly multiplicative, i.e., as maps $A\o A\o A\to K$,
\begin{eqnarray*}
(\eps\o \eps)\circ(m\o m)\circ(\id \o\cop\o\id )&=&\eps\circ
m\circ(m\o\id )\\ 
(\eps\o \eps)\circ(m\o m)\circ(\id \o\cop\op\o\id )&=&\eps\circ
m\circ(m\o\id ) 
\end{eqnarray*}
where $\cop\op:=\Sigma\circ\cop$ is opposite comultiplication,
\item $u$ is weakly comultiplicative, i.e., as maps $K\to A\o A\o A$,
\begin{eqnarray*}
(\id \o m\o\id )\circ(\cop\o\cop)\circ(u\o u)&=&(\cop\o\id )\circ\cop\circ
u\\
(\id \o m\op\o\id )\circ(\cop\o\cop)\circ(u\o
u)&=&(\cop\o\id )\circ\cop\circ u
\end{eqnarray*}
where $m\op:=m\circ\Sigma$ is opposite multiplication.
\end{enumerate}

Weak bialgebras reduce to ordinary bialgebras iff $\cop$ is unital.
Weak bialgebras have canonical subalgebras $A\L$ and $A\R$ that are spanned by
the right leg and left leg of $\cop(1)$, respectively. $A\L$ belongs to the
relative commutant of $A\R$ and there is a canonical antiisomorphism $A\L\to
A\R$. The subalgebras $A\L$ and $A\R$ are separable $K$-algebras.

Takeuchi's $\x_R$-bialgebras \cite{T} or, what is the same \cite{BM,Xu}
bialgebroids \cite{Lu} $A$ are defined over a ground ring $R$ which is not
supposed to be separable but plays the role of $A\L$ (or of $A\R$). Indeed
weak bialgebras are just the bialgebroids over separable base \cite{Sch2}. 

The weak bialgebras as well as the bialgebroids we speak about here are
finite dimensional over the ground field and finitely generated projective
over the base ring, respectively. Briefly saying they are finite quantum
groupoids.

\sect{Weak bialgebras versus bialgebroids}

A category of bialgebroids is introduced the arrows of which are called
bialgebroid maps. They intend to be special arrows in a possible larger
category. They have a special trend to point from bialgebras to bialgebroids
but not vice versa. The category of
maps of left/right bialgebroids will be denoted by $\Bia_l$, $\Bia_r$,
respectively. Using the forgetful functor from weak bialgebras to
bialgebroids we obtain weak morphisms of weak bialgebras as lifts of
bialgebroid maps. We comment on deformed versions of WBA's and on half
grouplike elements. 

\subsection{The category of bialgebroid maps}
Let $k$ be a commutative ring. All objects and maps below will be $k$-algebras
and $k$-algebra maps, respectively. Thus our base category is $k$-$\Alg$.

\subsubsection{The objects}
Following the original
definition \cite{T,Lu}, its reformulations in \cite{Xu,Sz}, and the terminology
of \cite{KSz} we say that
$\underline{A}=\bra A,R,s,t,\gamma,\pi\ket$ is a {\bf left
bialgebroid} if 
\begin{itemize}
\item $R\rarr{s}A\larr{t}R\op$ are $k$-algebra homomorphisms
such that $s(r')t(r)=t(r)s(r')$ for $r,r'\in R$. Then $A$ is made into an
$R$-$R$-bimodule by setting $r\cdot a\cdot r':=s(r)t(r')a$.
\item $\gamma\colon A\to A\o_R A$ and $\pi\colon A\to R$ are $R$-$R$-bimodule maps 
such that the triple $\bra A,\gamma,\pi\ket$ is a comonoid in the category
$_R\M_R$. 
\item $\gamma$ is a ring homomorphism into the Takeuchi $\x_R$ -product 
$A\x_R A$, i.e.,
\begin{eqnarray}
\gamma(a)(t(r)\o 1)&=&\gamma(a)(1\o s(r))\\
\gamma(a)\gamma(b)&=&\gamma(ab)\\
\gamma(1)&=&1\o 1
\end{eqnarray}
for all $a,b\in A$ and $r\in R$.
\item $\pi$ is compatible with the algebra structure, i.e.,
\begin{eqnarray}\label{Bia: eps}
\pi(as(\pi(b)))&=&\pi(ab)\ =\ \pi(at(\pi(b)))\,,\qquad a,b\in A\\
 \pi(1)&=&1_R\,.
\end{eqnarray}
\end{itemize}

\subsubsection{The arrows}

For two left bialgebroids $\underline{A}$ and $\underline{B}$ 
a pair $\bra \varphi,\omega\ket$ of algebra
homomorphisms $\varphi\colon A\to B$ and $\omega\colon R_A\to R_B$ is called a
\textbf{map of left bialgebroids} if 
\begin{equation}\label{dia 1-2}
\CD
A@>\varphi>>B\\
@A{s_A}AA@AAs_BA\\
R_A@>\omega>>R_B
\endCD
\qquad
\CD
A@>\varphi>>B\\
@A{t_A}AA@AAt_BA\\
R_A@>\omega>>R_B
\endCD
\end{equation}
are commutative diagrams in $k$-$\Alg$ and
\begin{equation}\label{dia 3}
\CD
A@>\varphi>>\Phi_\omega(B)\\
@V\pi_AVV@VV\pi_BV\\
R_A@>\omega>>R_B
\endCD
\end{equation}
\begin{equation} \label{dia 4}
\parbox[c]{4in}{
\begin{picture}(160,120)(-50,-20)
\put(0,80){$A$}
\put(10,84){\vector(1,0){110}}
\put(67,90){$\varphi$}
\put(122,80){$\Phi_\omega(B)$}
\put(3,75){\vector(0,-1){30}}
\put(-20,60){$\gamma_A$}
\put(-20,37){$A\o_{R_A}A$}
\put(7,32){\vector(2,-1){50}}
\put(5,13){$\varphi\o\varphi$}
\put(30,-5){$\Phi_\omega(B)\o_{R_A}\Phi_\omega(B)$}
\put(75,7){\vector(2,1){50}}
\put(100,13){$\tau^\omega_{B,B}$}
\put(100,37){$\Phi_\omega(B\o_{R_B}B)$}
\put(132,75){\vector(0,-1){30}}
\put(140,60){$\Phi_\omega(\gamma_B)$}
\end{picture}
}
\end{equation}
are commutative diagrams in $_{R_A}\M_{R_A}$.
Here $\bra\Phi^\omega,\tau^\omega,\omega\ket$ is the monoidal forgetful functor
$_{R_B}\M_{R_B}\to\,_{R_A}\M_{R_A}$ associated to the algebra homomorphism
$\omega$. That is to say for $\omega\colon R\to S$
\begin{eqnarray}
\Phi_\omega\colon& \,_S\M_S\to\,_R\M_R &\mbox{the forgetful functor}\\
\label{tau^omega}
\tau^\omega_{X,Y}\colon\Phi_\omega(X)\o_R\Phi_\omega(Y)&\to\Phi_\omega(X\o_S
Y)  &\mbox{the canonical}\\ 
\nonumber  &&\mbox{bimodule epimorphism}\\
\zeta_\omega\colon& \,_RR_R\to \,_RS_R &\mbox{$\omega$ as a
bimodule map.}
\end{eqnarray}

 Notice that $\omega$ is uniquely determined by $\varphi$ as
$\omega=\pi_B\circ\varphi\circ s_A$.

In order to see that the above properties of $\varphi$ are preserved under
composition of such maps, so we indeed have a category, one uses functoriality
of $\omega\mapsto \Phi^\omega$, i.e., that in fact the monoidal forgetful
functor is the arrow map of a functor $\Phi\colon \kAlg\op\to\MonCat$.
The object map of this functor is
\[
R\ \mapsto\ \bra\,_R\M_R,\o_R,\,_RR_R\,\ket
\]
Then functoriality means the identities
\begin{eqnarray*}
\Phi_\omega\circ\Phi_\sigma&=&\Phi_{\sigma\circ\omega}\\
(\Phi_\omega\circ\tau_\sigma)\bullet(\tau_\omega\circ(\Phi_\sigma\x\Phi_\sigma))
                        &=&\tau_{\sigma\circ\omega}\\
\Phi_\omega(\zeta_\sigma)\circ\zeta_\omega&=&\zeta_{\sigma\circ\omega}
\end{eqnarray*}
for $ R\rarr{\omega}S\rarr{\sigma}T$ in $\kAlg$.

In this way we have constructed a category $\Bia_l$ of left bialgebroids over
the base category $\kAlg$, i.e., the objects are left bialgebroids
\textit{in} $\kAlg$. In $\Bia_l$ there is no fixed base ring and there
are arrows between bialgebroids over different base rings. In particular there
are bialgebroid maps from ordinary bialgebras to bialgebroids.

In a similar way one defines the category $\Bia_r$ of right bialgebroids and
right bialgebroid maps the details of which we omit.

With the terminology "maps of bialgebroids" we indend to leave place for
more general arrows between bialgebroids. Certain bimodules with a coproduct,
so bimodule coalgebras, are natural candidates, they allow to
formulate Morita equivalence \cite{Sch3} when only the forgetful functor
$\M_A\to\,_R\M_R$ of the bialgebroid is considered as relevant. However, in the
Galois problem of non-commutative rings maps of bialgebroids do play a role 
as we shall see in Section \ref{sec: Gal}. In other words, $\Bia_l$ is large
enough to contain maps from group algebras or bialgebras to bialgebroids but
also small enough to contain only very restrictive isomorphisms.

\subsection{From weak bialgebras to bialgebroids}

Let $\underline{W}=\bra W,\cop,\eps\ket$ be a WBA over $K$. Its left and right
subalgebras are defined by
\begin{eqnarray}
L&=&\{\eps(1\1 w)1\2| w\in W\}\\
R&=&\{1\1\eps(w 1\2)| w\in W\}
\end{eqnarray}
and are the images of the maps
\begin{eqnarray}
\PL\colon W\to L\,,\qquad&& w\mapsto \eps(1\1 w)1\2\\
\PR\colon W\to R\,,\qquad&& w\mapsto 1\1\eps(w 1\2)\ .
\end{eqnarray}
For the basic properties of these maps see \cite{BNSz}.

Now we introduce data for bialgebroids as follows. Let
\begin{eqnarray}
s\L\colon &&L\to W\\
s\R\colon &&R\to W
\end{eqnarray}
just the inclusion maps, hence algebra homomorphisms. Let
\begin{eqnarray}
t\L\colon &&L\op\to W\,,\qquad l\mapsto 1\1\eps(1\2 l)\\
t\R\colon &&R\op\to W\,,\qquad r\mapsto \eps(r 1\1)1\2
\end{eqnarray}
which are also algebra maps (if antipode exists they are the restrictions of
$S^{-1}$). The ranges of $s\L$ and $t\L$ are $L$ and $R$, so they commute.
Similarly for $s\R$ and $t\R$. This allows us to introduce bimodule structures
$_LW_L$ and $_RW_R$, respectively, via the formulae
\begin{eqnarray}
l\cdot w\cdot l'&:=&s\L(l)t\L(l')w\,,\qquad l,l'\in L,\ w\in W\\
r\cdot w\cdot r'&:=&w t\R(r)s\R(r')\,,\qquad r,r'\in R,\ w\in W\ .
\end{eqnarray}
Finally, let
\begin{eqnarray}
\tau\L \colon && W\o W\to W\o_L W\\
\tau\R \colon && W\o W\to W\o_R W
\end{eqnarray}
be the canonical epimorphisms associated to the units $K\to L$,
$K\to R$, repectively. 

\begin{lemma}
Let $\underline{W}=\bra W,\cop,\eps\ket$ be a WBA$/K$. Then the maps
\begin{eqnarray}
\gamma\L:=\tau\L\circ\cop\colon &&W\to W\o_L W\\
\gamma\R:=\tau\R\circ\cop\colon &&W\to W\o_R W
\end{eqnarray}
are such that 
\begin{eqnarray}
\beta_l(\underline{W})&:=&\bra W,L,s\L,t\L,\gamma\L,\PL\ket\ \mbox{is a left
bialgebroid}\\
\beta_r(\underline{W})&:=&\bra W,R,s\R,t\R,\gamma\R,\PR\ket\ \mbox{is a right
bialgebroid.}
\end{eqnarray}
\end{lemma}
{\it Proof}:\,\,
First notice that $\pi\L\circ s\L=\id_L$ and $\pi\L\circ t\L=\id_L$ by
\cite[(2.3a), (2.25a)]{BNSz}. Also using \cite[Lemma 2.5]{BNSz} 
\begin{eqnarray*}
\pi\L(l\cdot w\cdot l')&=&
\pi\L(lt\L(l')w)=l\pi\L(t\L(l')\pi\L(w))=l\pi\L(w)\pi\L(t\L(l'))\\ 
&=&l\pi\L(w)l'
\end{eqnarray*}
therefore $\pi\L$ is an $L$-$L$ bimodule map. Now $\bra
\,_LW_L,\gamma\L,\pi\L\ket$ is a comonoid in $_L\M_L$ because $\tau\L$
is natural and $\pi\L(w\1)w\2=w$ and $t\L\circ\pi\L(w\2)w\1=w$ are general WBA
identities. Coassociativity follows using that $\tau\L$ extends to a 
natural transformation (\ref{tau^omega}), namely, $\tau\L=\tau^{u_L}_{W,W}$
where $u_L\colon K\to L$ is the unit of the $K$-algebra $L$, and the latter
satisfies the hexagon diagram of a (lax) monoidal functor. In order to show
that the image of $\gamma\L$ is in $W\x_L W$ it suffices to refer to the old
WBA identity $1\1\o 1\2 l=1\1 t\L(l)\o 1\2$ (cf. \cite[(2.31a)]{BNSz}).
Multiplicativity of $\gamma\L$ then follows from multiplicativity of $\cop$.
It remains to show the counit properties (\ref{Bia: eps}) but they are just the
identities \cite[(2.5a),(2.25a)]{BNSz}. Passing to the opposite-coopposite
WBA one obtains the statement for $\beta_r$.
\hfill{\it q.e.d.}

The $\beta_l$ and $\beta_r$ defined by the Lemma are expected to be the object
maps of two functors
\begin{equation}
\begin{CD}
\Bia_l@<\beta_l<<WBA@>\beta_r>>\Bia_r
\end{CD}
\end{equation}
the arrows of WBA, however, will be discussed later.

\subsection{From bialgebroids over separable base to weak bialgebras}

Since the left/right subalgebras of a WBA are always separable, we start from
a left bialgebroid $\underline{B}=\bra B,R,s,t,\gamma,\pi\ket$ in which $R$ is
a separable $K$-algebra. That is to say there exists an element
$e=\sum_i e_i\o e^i\in R\o_K R$ such that
\begin{equation}
\sum_i e_ie^i=1_R\quad\mbox{and}\quad \sum_i re_i\o e^i=\sum_ie_i\o e^ir\quad
r\in R\,.
\end{equation}
Such a separability idempotent provides a splitting map for the canonical
epimorphism $\tau$, namely
\begin{equation}\label{sigma}
\sigma\colon B\o_R B\to B\o_K B\,,\quad b\o b'\mapsto\sum_i b\cdot e_i\o
e^i\cdot b'
\end{equation}
This formula is the same for right bialgebroids. For left bialgebroids we can
write also
\begin{equation}
\sigma(b\o_R b')=\sum_i t(e_i)b\o s(e^i)b'\,.
\end{equation}

But there is more than separability of $R$ in a WBA. There is also a
separabity \textit{structure} for $R$. For any weak bialgebra with left
subalgebra $R$ the restriction of the counit  $\psi:=\eps|_R$ is a
nondegenerate functional of index one. This means that $\psi$ distinguishes a
special separability idempotent, namely $e=S(1\1)\o 1\2$ which is the
quasibasis of $\psi$. Comparing this with the above expression for $\sigma$
one recognizes that $\sigma$ is multiplication from the left by $\cop(1)$ on
any element from the inverse image $\tau^{-1}(\{b\o_R b'\})$. 

\begin{lemma} \label{lem: biatoWBA}
Given a pair $\bra\underline{B},\psi\ket$, where $\underline{B}$ is a left or
right bialgebroid over $R$ and $\psi\colon R\to K$ is a nondegenerate
functional of index 1, define
\begin{eqnarray}
\cop&:=&\sigma\circ\gamma\colon B\to B\o_K B\\
\eps&:=&\psi\circ\pi\colon B\to K
\end{eqnarray}
where $\sigma$ is the splitting map of the canonical epimorphism $B\o_K B\to
B\o_R B$ that is associated to the quasibasis $e$ of $\psi$ as in
(\ref{sigma}). Then the triple $\bra B,\cop,\eps\ket$ is a WBA over $K$.
\end{lemma}
{\it Proof}:\,\,
Mutatis mutandis, the proof has already been given in \cite[Proposition
9.4]{KSz} and in \cite[Theorem 5.5]{Sch2}. 
\hfill{\it q.e.d.}

The above Lemma characterizes the fibres of the functor $\beta_l$ in the
following sense. The WBA's $W$ with a fixed underlying (let's say left)
bialgebroid $\beta_l(W)=B$ are in one-to-one correspondence with separability
structures $\bra R,\psi,e\ket$ on $R$.

\subsection{Strict and weak morphisms of weak bialgebras}

\begin{definition} \label{def: str/weak mor}
Let $\bra W,\cop,\eps\ket$ and $\bra W',\cop',\eps'\ket$ be weak bialgebras
over $K$. Then a $K$-linear map $f\colon W\to W'$ is called a
\begin{itemize}
\item {\bf strict morphism} of weak bialgebras if $f$ is an algebra map and a
coalgebra map;
\item {\bf weak left morphism} of weak bialgebras if
$f\colon\beta_l(W)\to\beta_l(W')$ is a map of the underlying left bialgebroids;
\item {\bf weak right morphism} of weak bialgebras if
$f\colon\beta_r(W)\to\beta_r(W')$ is a map of the underlying right
bialgebroids.
\end{itemize}
\end{definition}
A strict morphism $f$ not only preserves the left and right subalgebras,
$f(L)\subset L'$ and $f(R)\subset R'$, but, because $(f\o
f)\circ\cop(1)=\cop'(1')$, establishes isomorphisms $L\iso L'$ and $R\iso R'$.
Therefore strict morphisms exist between two WBA's only if they have isomorphic
left, resp. right subalgebras. This is definitely too strong since the
original philosophy of \cite{BSz} was to "blow up" Hopf algebras in order they
could afford non-integral categorical dimensions, but the amount of the
blowing up, i.e., the size of the left/right subalgebras should be considered
as a gauge degree of freedom. Using weak morphisms we pursue this
idea to some extent. 

Since weak morphisms are justs lifts of the rather involved bialgebroid maps
into the WBA framework, they are useful only if they can be recognized
directly without reference to bialgebroids. Therefore we make the

\begin{proposition} \label{pro: wlmor}
For weak bialgebras $\underline{W}$ and $\underline{W'}$ a $K$-linear map
$f\colon W\to W'$ is a weak left morphism iff 
\begin{enumerate}                           
\item $f$ is a $K$-algebra map, \label{wlm 1}
\item $f(R)\subset R'$, \label{wlm 2}
\item ${\pi'}\L\circ f=f\circ\pi\L$, \label{wlm 3}
\item and $\cop'(1')\cdot(f\o f)(\cop(w))=\cop'(f(w))$, $w\in W$.\label{wlm 4}
\end{enumerate}
It is a weak right morphism iff
\begin{enumerate}                           
\item $f$ is a $K$-algebra map,
\item $f(L)\subset L'$,
\item ${\pi'}\R\circ f=f\circ\pi\R$,
\item and $(f\o f)(\cop(w))\cdot\cop'(1')=\cop'(f(w))$, $w\in W$.
\end{enumerate}
\end{proposition}
{\it Proof}:\,\,
It suffices to prove the statement for left morphisms. Since $s\L$ is just the
injection $L\subset W$ and so is ${s'}\L$, the first diagram in (\ref{dia
1-2})  is equivalent to $f(L)\subset L'$, which in turn is a consequence of
(\ref{dia 3}) which is nothing but condition \ref{wlm 3} above. Having
condition \ref{wlm 3} anyway the second diagram of (\ref{dia 1-2}) is
equivalent to the condition \ref{wlm 2} because if \ref{wlm 2} holds
then $f\circ t\L={t'}\L\circ{\pi'}\L\circ f\circ t\L={t'}\L\circ
f\circ\pi\L\circ t\L={t'}\L\circ f$ and backwards is obvious. This proves that
the three diagrams of (\ref{dia 1-2}) and (\ref{dia 3}) are equivalent to
\ref{wlm 1}, \ref{wlm 2} and \ref{wlm 3}. Assuming this we can equip $W'$ with
$L$-$L$ bimodule structure $\hat W'=\Phi_\omega(\,_{L'}W'_{L'})$ and lift $f$
to an $L$-$L$ bimodule map $\hat f$. Then (\ref{dia 4}) takes the form
\begin{equation}\label{eq: dia 4 again}
\tau^\omega\circ (\hat f\o_L\hat f)\circ\gamma=\gamma'\circ\hat f\,.
\end{equation}
Let $\tau'$
be the canonical epi for $W'$ and $\sigma'$ be its splitting map associated to
the quasibasis $e'$ of $\eps'|_{L'}$. Let $\hat\tau$ be the canonical epi for
$\hat W'$ and $\hat\sigma$ its splitting map that is associated to $e$, or
to $(f\o f)(e)$ in some (bad) sense. Then we have $\tau^\omega\circ\hat\tau=
\tau'$ and using the observation we made just before Lemma \ref{lem: biatoWBA}
we can write for all $w'_1\o w'_2\in \hat W'\o_L \hat W'$ that
\begin{eqnarray*}
\sigma'\circ\tau^\omega(w'_1\o w'_2)&=&\sigma'\circ\tau^\omega\circ
\hat\tau\circ\hat\sigma(w'_1\o w'_2)\\
&=&1'\1f(1\1)w'_1\o 1'\2f(1\2)w'_2\,.
\end{eqnarray*}
Now acting by $\sigma'$ on (\ref{eq: dia 4 again}) we obtain \ref{wlm 4} and
acting by $\tau'$ on \ref{wlm 4} we obtain (\ref{eq: dia 4 again}).
\hfill{\it q.e.d.}

For weak Hopf algebras one defines weak left/right morphisms as those of its
underlying weak bialgebra, disregarding whether they preserve antipodes or
not.

The category of bialgebras, as well as the category of Hopf algebras,
are full subcategories in each one of $\Bia_l$, $WBA$ and $\Bia_r$.

\begin{example}
Let $H$ be a Hopf algebra over $K$. Define its blowing up as the algebra
$W:=H\o M_n(K)$ with comultiplication $\cop(h\o e_{ij}):=(h\1\o
e_{ij})\o(h\2\o e_{ij})$. Then $W$ becomes a weak Hopf algebra. Its left and
right subalgebras coincide and equal to the diagonal matrices with entries
from $K$. The diagonal embedding of $H$, $f\colon H\to W$, $f(h)=h\o I_n$, is
clearly an algebra map. It is not a coalgebra map however, but we have
\begin{eqnarray*}
\cop(f(h))&=&\sum_i (h\1\o e_{ii})\o (h\2\o e_{ii})\\
          &=&\cop(1_W)(f\o f)(\cop_H(h))=(f\o f)(\cop_H(h))\cop(1_W)\\
\pi\L(f(h))&=&f(1_H)\eps_H(h)=\pi\R(f(h))\,.
\end{eqnarray*}
Now using the above Proposition it is plain that $f$ is weak left and right
morphism of weak bialgebras and there is no strict morphism from $H$ to $W$
unless $n=1$.
\end{example}

\subsection{Weak automorphisms, twists and half grouplike elements}

Let $f\colon W\to W'$ be a weak left morphism of WBA's. Then
$\eps'(f(w))=\eps'({\pi'}\L(f(w)))=\eps'(f(\pi\L(w))=\eps(1\1 w)\eps'(f(1\2))$
therefore
\begin{equation}
\eps'(f(w))=\eps(uw)\,,\quad w\in W,\qquad\mbox{where } u=\eps'(f(1\1))1\2\in L
\end{equation}
where, in order to get $u\in L$, we also made the $R\to L$ transformation
$t\L(1\1\eps'(f(1\2)))=\eps'(f(\pi\L(1\1)))1\2=u$.
Especially we have $\eps'(f(l))=\eps(ul)$ for $l\in L$. So assuming that
$f|_L$ is an isomorphism onto $L'$ we have $u$ as a Radon-Nikodym derivative
of a nondegenerate functional w.r.t another, hence invertible. Comparing their
quasibases we obtain the equality
\begin{equation}
\omega^{-1}({\pi'}\L(1'\1))\o \omega^{-1}(1'\2)=\pi\L(1\1)\o u^{-1}1\2
\end{equation}
as elements of $L\o L$. Applying $f\o f$ and using that $\pi\L$ restricts to
an isomorphism $R\to L$ (the would-be antipode), we get
\begin{equation}
1'\1\o 1'\2=f(1\1)\o f(u^{-1})f(1\2)\,.
\end{equation}
Inserting this result to the 4th property of weak left morphisms in
Proposition \ref{pro: wlmor} one immediately arrives to the
\begin{lemma} \label{lem: defo}
Let $f\colon W\to W'$ be a weak left morphism of WBA's such that its
restriction $\omega\colon L\to L'$ is an isomorphism. Then there is an
invertible element $u\in L$ such that 
\begin{eqnarray}
\cop'(f(w))&=&(f\o f)((1\o u^{-1})\cop(w))\\
\eps'(f(w))&=&\eps(uw)
\end{eqnarray}
for all $w\in W$.
\end{lemma}

This result holds in particular if $f$ is an isomorphism. As a matter of fact
property 3 in Proposition \ref{pro: wlmor} is invariant under changing $f$
to $f^{-1}$. For completeness we remark that the forgetful functor $\Bia_l\to
\kAlg$ reflects isomorphisms. That is to say, if a weak left morphism is
invertible as an algebra map then its inverse is a weak left morphism. 
So the Lemma holds for $f=\id_W$. In this case $u$ describes a deformation in
the sense of \cite[Remark 3.7]{Nik}. Such (left) deformed WBA's have identical
underlying left bialgebroids, so deformations should be interpreted as weak
left automorphisms. Although the deformation changes the Nakayama
automorphism of the counit, there may be no deformation at all which
produces a tracial $\eps'$, unless the base $L$ possesses a nondegenerate trace
of index 1. For example if $L$ is split semisimple then the only such trace is
the regular trace. Since the Radon-Nikodym derivative of the regular trace
w.r.t. $\eps|_L$ is $1\2S(1\1)$, tracial deformation exists iff $1\2S(1\1)$ is
invertible. In characteristic zero this is always the case, otherwise there
are counter examples \cite{V}.

Now consider inner weak left automorphisms associated to left grouplike
elements. For a bialgebroid $\bra B,L,s,t,\gamma,\pi\ket$ an element $g\in B$
is grouplike if $g$ is invertible and $\gamma(g)=g\o_L g$. For a WBA $W$ define
\begin{eqnarray}
\G\L(W)&:=&\{g\in W| g\mbox{ is grouplike in }\beta_l(W)\}\\
\G\R(W)&:=&\{g\in W| g\mbox{ is grouplike in }\beta_r(W)\}
\end{eqnarray}
the sets of left/right grouplike elements.
This is of course equivalent to saying e.g. that $g$ is left grouplike if it
is invertible and $\cop(g)=\cop(1)(g\o g)$. In the next computations we assume
that $g,h\in\G\L(W)$ and $w\in W$ is arbitrary.
\begin{eqnarray}
&&\pi\L(g)=\eps(1\1 g)1\2=\eps(g\1)g\2 g^{-1}=gg^{-1}=1\\
&&\cop(g^{-1})=\cop(g^{-1})\cop(1)=\cop(g^{-1})\cop(g)(g^{-1}\o g^{-1})\nonumber\\
&&\quad=\cop(1)((g^{-1}\o g^{-1})\\
&&\cop(gh)=\cop(g)\cop(1)(h\o h)=\cop(g)(h\o h)=\cop(1)(gh\o gh)\\
&&\pi\L(gwg^{-1})=\pi\L(gw\pi\L(g^{-1}))=\pi\L(gw)=\eps(1\1
g\pi\L(w))1\2\nonumber\\
&&\quad=\eps(g\1\pi\L(w))g\2 g^{-1}=g\pi\L(w)g^{-1}\\
&&t\L\circ\pi\L(gwg^{-1})=1\1\eps(1\2 gwg^{-1})=g\1
g^{-1}\eps(g\2\pi\L(w))\nonumber\\
&&\quad=g(t\L\circ\pi\L(w))g^{-1}\\
&&\cop(1)(gw\1 g^{-1}\o gw\2 g^{-1})=g\1 w\1 g^{-1}\o g\2 w\2 g^{-1}\nonumber\\
&&\quad=g\1 w\1 g^{-1}\1\o g\2 w\2 g^{-1}\2=\cop(gwg^{-1})
\end{eqnarray}
This shows that $\G\L$ is a group and for all $g\in\G\L$ the inner
automorphism $w\mapsto gwg^{-1}$ is a weak left automorphism $W\to W$.
It does not leave the counit invariant but
\begin{equation}
\eps(gwg^{-1})=\eps(uw)\quad w\in W\,,\qquad \mbox{where } u=\eps(g1\1)1\2
\end{equation}
implying for their quasibasis the relation
\[
g^{-1}\pi\L(1\1)g\o g^{-1}1\2 g=\pi\L(1\1)\o u^{-1}1\2\,.
\]
Applying $t\L\o\id$ and using that $\Ad_{g^{-1}}$ commutes with $\pi\L$ and
$t\L\circ\pi\L$ one obtains 
\[
g^{-1}1\1 g\o g^{-1} 1\2 g=1\2\o u^{-1} 1\2
\]
therefore
\begin{equation}
\cop(g)=(g\o gu^{-1})\cop(1)=(g\pi\R(u^{-1})\o g)\cop(1)\ .
\end{equation}
Since $\pi\R$ on $L$ is an algebra antiisomorphism and $\pi\R(u)=1\1\eps(g
S^{-1}(1\2))$ $=\pi\R(g)$, the above definition of $\G\L$ is equivalent to the
one given by \break Vecserny\'es \cite{V}.

\sect{Galois quantum groupoids} \label{sec: Gal}

In this section we argue that the balanced depth 2 extensions \cite{KSz} of
rings or $k$-algebras are the proper analogues of the Galois extensions of
fields (i.e., normal and separable field extensions) because they have finite
quantum automorphism groups (cf. \cite{Pareigis}) with invariant subalgebra
just $N$ and which are characterized by a universal property, hence unique. The
role of finite groups are played by bialgebroids, i.e., $\x_R$-bialgebras,
that are finitely generated projective over their base as a left and as a
right module. They will be called finite bialgebroids. They are presumably
also Hopf algebroids but the antipode raises several questions, so we skip
their discussion altogether. The difference between Galois bialgebroid and
Galois WBA will be found in the difference between depth 2 Frobenius
extensions and depth 2 extensions with a Frobenius structure.

\subsection{Quantum automorphisms}

Recall the definition of left module algebroids
over a left bialgebroid $A$ in \cite{KSz}. They are the monoids in the
category of left $A$-modules. 

\begin{definition} \label{def: Gal}
Let $N\to M$ be an extension of $k$-algebras. Define the category $\Aut(M/N)$
as follows. Its objects are the pairs $\bra B, \alpha_B\ket$ where $B$ is a
finite left bialgebroid in $k$-$\Alg$ and $\alpha_B\colon B\o M\to M$ is a
left $B$-module algebroid action such that $N$ is contained in the invariant
subalgebra $M^B$. The arrows from $\bra B,\alpha_B\ket $ to $\bra
C,\alpha_C\ket$ are the maps $f\colon B\to C$ of left bialgebroids such that
\[
\alpha_C\circ(f\o \id_M)=\alpha_B\ .
\]
A terminal object in $\Aut(M/N)$ is called a \textbf{universal action} on the
extension $M/N$. The bialgebroid $A$, unique up to isomorphism in $\Bia_l$, in
a universal action is called the \textbf{Galois quantum groupoid} of the
extension $M/N$ and is denoted by $\Gal(M/N)$. If the invariant subalgebra of
the universal  action is equal to $N$ the extension $N\to M$ is called a
\textbf{Galois extension}. 
\end{definition}

If $F\subset E$ is a field extension then every finite group action on $E$
which leaves $F$ pointwise fixed factors uniquely through the Galois group
$\Gal(E/F)$. This trivial fact is generalized by the above definition. 
Also, if the fixed points of $\Gal(E/F)$ coincide with the elements of $F$ the
extension is normal and separable, i.e., Galois, by Artin's Theorem.

Note that the real beauty of universal monoids of \cite{Pareigis} has not been
used in the above definition. One could consider much more general arrows
$\alpha\colon B\o M\to M$ than just actions. 

If $H$ is a Hopf algebra, f.g.p. over $k$ then Kreimer and Takeuchi defines an
$H$-Galois extension to be a ring extension $N\subset M$ such that 
\begin{itemize}
\item there is a left $H$-module algebra action $\alpha\colon H\o M\to M$,
\item $N=M^H$, the invariant subring,
\item $M_N$ is finitely generated projective and 
\item the map 
\begin{equation} 
\phi\colon M\o H\to \End M_N\,,\quad m\o h\mapsto \{m'\mapsto m\alpha(h\o m')\}
\end{equation}
is an isomorphism. 
\end{itemize}
(More precisely, this is a reformulation by Ulbrich \cite{Ulbrich}.)

Hopf-Galois extensions in this sense, however, do not have the universal
property with respect to the category of Hopf algebras. As it was pointed out
by Greither and Pareigis in \cite{Gr-Pa} there are separable field extensions
which are $H$-Galois for two different Hopf algebras. We come back to this
example in the last section.

\subsection{Universal bialgebroid actions}

The advantage of using bialgebroids is that there is a very general class of
ring extensions for which a universal bialgebroid action exists. These are the
depth 2 extensions $N\subset M$ for which $M_N$ is balanced. They include all
ring extensions that are $H$-Galois for some Hopf algebra, as it was shown by
Kadison recently \cite{K2}, but many more. The universal bialgebroid of a
depth 2 balanced ring extension is a canonical structure on the endomorphism
ring $A=\End\,_NM_N$ and it has been constructed in \cite{KSz} although its
universality was not formulated there. Below we shall give a proof for the
special case of separable centralizer which leads us to weak bialgebra actions
as follows. 

If $W$ is a weak bialgebra over $K$ and $B=\beta_l(W)$ its underlying left
bialgebroid then the category of left $W$-modules and the category of left
$B$-modules are monoidally equivalent \cite{Sch2}, Proposition 5.3, in fact
isomorphic. Therefore these categories have the same monoids. Therefore a
module algebra over $W$ is the same as a module algebroid over $B$. This lends
to a WBA action $\alpha\colon B\o M\to M$ the name {\bf weak Galois action} if
it has $N$ as its invariant subalgebra and if it has the universal property
w.r.t. weak left morphisms of WBA's.

\begin{theorem} \label{thm: unibia}
Let $K$ be a field and $N\subset M$ a $K$ algebra extension such that 
\begin{itemize}
\item $N\subset M$ is of depth 2,
\item $M_N$ is balanced,
\item $R:=C_M(N)$ is a separable $K$-algebra.
\end{itemize}
Then the bialgebroid $A=\End\,_NM_N$ constructed in \cite{KSz} and acting
on $M$ in the natural way is the Galois bialgebroid $\Gal(M/N)$. That is to
say, any weak bialgebra $\bra A,\cop_A,\eps_A\ket$ with underlying left
bialgebroid $A$ has the following universal property. If $\alpha_W\colon W\o
M\to M$ is a left module algebra action of a WBA $W$ such that
$M^W\supset N$ then there is a unique weak left morphism $f\colon W\to A$ of
weak bialgebras such that $\alpha\circ(f\o\id_M)=\alpha_W$.
\end{theorem}

{\it Proof}:\,\,
That the invariant subalgebra is $N$ was shown in \cite{KSz}. To prove the
universal property notice that for any weak bialgebra action on $M$
\begin{equation}\label{eq: nmn}
w\lact(nmn')=n(w\lact m)n'\,,\quad w\in W,\ n,n'\in M^W,\ m\in M\,,
\end{equation}
in particular for all $n,n'\in N$. Thus there is a unique algebra map $f\colon
W\to A$ such that 
\[
f(w)(m)=w\lact m\,.
\]
Already this implies uniqueness so we are left to show that $f$ is a weak left
morphism. We will use the criteria given in Proposition \ref{pro: wlmor}.
At first compute the action of an $l\in W\L$.
\[
l\lact m=(l\1\lact 1_M)(l\2\lact m)=1\lact((l\lact 1_M)m)=(l\lact 1_M)m
\]
Since $(l\lact 1_M)n=l\lact n=n(l\lact 1_M)$ for $n\in N$, $W\L\lact
1_M\subset C_M(N)=R$. Hence $f(W\L)\subset\lambda_M(R)=A\L$, where
$\lambda_M(m)$ denotes left multiplication by $m$ on $M$. For $r\in W\R$ we
have
\[
r\lact m=(r\1\lact m)(r\2\lact 1_M)=1\lact(m(r\lact 1_M))=m(r\lact 1_M)
=m(\pi_W\L(r)\lact 1_M)
\]
Since $A\R=\rho_M(R)$ where $\rho_M(m)$ is right multiplication by $m$ on $M$,
we obtain $f(W\R)\subset\rho_M(R)=A\R$.
Next recall that the counit of the left bialgebroid $A$ is
$\pi_A\L(a)=\lambda_M(a(1_M))$. Therefore
\[
\pi_A\L (f(w))=\lambda_M(f(w)(1_M))=\lambda_M(w\lact
1_M)=\lambda_M(\pi_W\L\lact 1_M)=f(\pi_W\L(w))
\]
Turning to the last condition of Proposition \ref{pro: wlmor} we recall
\cite{KSz}, Prop. 3.9 stating that
\[
A\o_R A\iso \Hom_{N-N}(M\o_N M,M),\quad (a\o a')(m\o m')=a(m)a'(m')\,.
\]
Now consider the composite $K$-linear maps
\[
\begin{CD}
W@>\cop>>W\o W@>f\o f>>A\o A@>\tau>> A\o_R A\\
W@>f>>A@>\gamma>>A\o_R A
\end{CD}
\]
which are equal because their images act on $M\o_N M$ in the same way due to
the module algebra property. So, composing them with the splitting map
$\sigma$ associated to $\cop_A(1_A)$ and using $\sigma\circ\tau=\lambda_{A\o
A}(\cop_A(1_A))$ we obtain
\[
\cop_A(1_A)\cdot (f\o f)(\cop_W(w))=\cop_A(f(w))\,,\quad w\in W\,.
\]
\hfill{\it q.e.d.}

\subsection{Universal weak Hopf algebra actions}

If we add to the conditions of Theorem \ref{thm: unibia} that $N\subset
M$ is a Frobenius extension then it already implies that the WBA lift of its
Galois bialgebroid is a WHA \cite[Section 9]{KSz}. This does not make the WHA
unique, only up to weak left isomorphisms. This freedom of the WHA is
precisely the freedom of choosing a Frobenius functional $\psi\colon R\to K$
of index 1. Therefore it is natural to associate WHA actions to Frobenius
structures $\bra N\to M, \phi\colon \,_NM_N\to\,_NN_N\ket$ rather than to just
extensions $N\to M$. 

If $\phi\colon\,_NM_N\to\,_NN_N$ is a Frobenius map, i.e., a bimodule map with
quasibasis $\sum_i m_i\o m^i\in M\o_N M$, then
its restriction to the centralizer $\phi|_R$ maps $R$ into the center $Z$ of
$N$. Since $R$ is not only part of $M$ but belongs to the bialgebroid $A$ as
well, it is very natural to build $\phi|_R$ into the data of the WHA as
$\eps|_R$. Strictly speaking, this is possible only if the center of $N$ is
trivial. There is a tiny point here about the
restriction. While in case of finite index $C^*$-algebra extensions one
considers faithful conditional expectations $\phi$ which have faithful
restrictions to the finite dimensional $R$, therefore $\phi|_R$ is a Frobenius
map with invertible index, this is not automatic for general Frobenius algebra
extensions. 

\begin{theorem} \label{thm: uniWHA}
Let $N\subset M$ be a depth 2 Frobenius extension of $K$-algebras with
centralizer $R$ a separable $K$-algebra and with $\Center N=K$. Assume
$\phi\colon M\to N$ is a Frobenius map with its restriction $\phi|_R$ being an
index 1 Frobenius map. Then there exists a unique weak Hopf algebra $A$ and a
left module algebra action of $A$ on $M$ which satisfies the universal
property of Theorem \ref{thm: unibia} and such that $\phi|_R=\eps|_R$. 
\end{theorem}
{\it Proof}:\,\,
The antipode of a WHA is unique therefore uniqueness of $A$ follows if we show
that its WBA structure is unique. The latter is uniquely determined by its
underlying left bialgebroid $\beta_l(A)$ and by the restriction of its counit,
$\eps|_R$. The former is uniquely determined by the universal property as
$\beta_l(A)=\Gal(M/N)$ by Theorem \ref{thm: unibia} and the latter by the
requirement $\eps|_R=\phi|_R$. This proves uniqueness. The existence part is
an easy application of Theorem 9.5 of \cite{KSz}. 
\hfill{\it q.e.d.}

The question arises how to interpret $\phi$ if only its restriction to the
centralizer matters. Since $\phi$ is an $N$-$N$ bimodule map, it belongs to
$A$ as a nondegenerate left integral. Thus in fact the data of the Theorem
determine a measurable quantum groupoid, i.e., a WHA with a distinguished
nondegenerate integral. 

Generalizations to $\Center N=Z$ a separable $K$-algebra is possible. It
requires to use a slight generalization of the notion of a WHA. It requires
WHA's not in $\kAlg$ but in $_Z\M_Z$, cf. \cite[Proposition 1.6]{Sz}.

In addition to the assumptions of Theorem \ref{thm: uniWHA} let us assume that
$\phi(1_M)$ is invertible or only assume that $M/N$ is split. Then $M_N$ is
balanced therefore $M/N$ is a Galois extension in the sense of Definition
\ref{def: Gal}.

\sect{Separable field extensions are weak Hopf Galois}

Let $E|K$ be a separable field extension. Then the following results are
standard.
\begin{enumerate}
\item
$E_K$ is finite dimensional, $\dim E_K=n<\infty$.
\item Let $\tau\colon E\to K$ be the trace associated to the regular
$E$-module $_EE$. Then there exists $x_i,y_i\in E$, $i=1,\dots n$ such that
\begin{eqnarray}
\sum_i\,x_i\tau(y_ix)&=&x\,,\qquad \forall x\in E\\
\sum_i\ x_iy_i&=&1\,.
\end{eqnarray}
\item
The above set $\{x_i,y_i\}$, called a quasibasis for $\tau$, satisfies
\[
\sum_i\,xx_i\o y_i\ =\ \sum_i\,x_i\o y_ix\,,\qquad \forall x\in E\,.
\]
\item
There exists a $\xi\in E$ which generates $E$ as an $K$-algebra, i.e.,
$E=K(\xi)$. (Primitive Element Theorem)
\item
The non-zero $K$-algebra endomorphisms of $E$ are automorphisms. \break Their
group $G$ forms an $K$-linearly independent set in the $K$-algebra of
$K$-linear endomorphisms $\End E_K$ of the $K$-module $E_K$. The
$G$-invariants $F=E^G$ form a subfield of $E$ and $E|F$ is (classically)
Galois with Galois group $G$. Hence $|G|=\dim_FE=n/m$ where $m=\dim_KF$.
\end{enumerate}

\subsection{The universal weak Hopf algebra of $E/K$}

Define $A$ as the $K$-algebra $\End E_K$ and its weak Hopf algebra structure by
\begin{eqnarray}\label{WHA A}
\cop_A(a)&=&\sum_i\sum_j\ x_i\tau(x_j\under)\ \o\ y_ja(y_j\under)\,,\\
\eps_A(a)&=&\tau(a(1))\,,\\
S_A(a)&=&\sum_i x_i\tau(a(y_i)\under)\,.
\end{eqnarray}
The WHA $A$ is a very special one:
\begin{enumerate}
\item The left and right subalgebras coincide with $E$. As a matter of fact,
identifying $E$ with the subalgebra of $A$ of (left) multiplications on $E$
\[
\pi^L(a)=a(1)\,,\qquad\pi^R(a)=\sum_ix_i\tau(a(y_i))
\]
thus $A^L=A^R=E$.
\item The antipode is involutive, $S_A^2=\id_A$. What is more, it is
transposition w.r.t. the nondegenerate bilinear form on $E\o_K E$ given by
$\tau$, i.e., 
\[
\tau(xa(y))\ =\ \tau(S_A(a)(x)y)\,,\qquad x,y\in E\,.
\]
\item $A$ is cocommutative, $a\1\o a\2=a\2\o a\1$ as elements of $A\o_K A$,
holds for all $a\in A$ as a consequence of commuativity of $E$ on the one hand
and of the existence of an isomorphism $A\o_K A\cong \End_K(E\o_KE)$, i.e.,
finite dimensionality of $E_K$ on the other hand.
\item Left integrals in $A$ are the endomorphisms $l$ of $E_K$ such that
$l(E)\subset K$. Their general form is $l=\tau(r\under)$ where $r\in E$.
Normalized left integrals thus exist ($\Rightarrow\ A$ is a separable
$K$-algebra) and the invariant subalgebra of $E$ is $K$.
\item $\tau$ is a 2-sided nondegenerate integral. If $n$ is invertible 
in $K$, especially in characteristic 0, then $\tau/n$ is a Haar integral in $A$.
\item The left grouplike elements of $A$ are precisely the algebra
automorphisms of $E_K$. Thus the number $|\G\L|$ of left grouplike elements
is a divisor of $n$ and is equal to $n$ precisely if $E/K$ is
classically Galois. In the latter case $A$ is a crossed product of $E$ with
the group algebra $K\G\L$.  
\item The left $A$-module algebra $A\L$, i.e., the trivial
$A$-module, coincides with $E$ with its canonical left $A$-module structure.
That is to say, 
\begin{eqnarray*}
a(x)&=&ax(1)=a\1x(1\eps_A(a\2))=a\1xS_A(a\2)(1)=\pi^L(ax)(1)\\
&=&\pi^L(ax)
\end{eqnarray*}
for all $a\in A$ and $x\in E$.
\item The smash product $E\rcross A$ is isomorphic to $A$ as $K$-algebras via
the canonical map $x\rcross a\mapsto\{y\mapsto xa(y)\}$. 
\end{enumerate}

\subsection{Weak Hopf Galois extensions}

As an immediate generalization of Hopf-Galois extensions one can make the
\begin{definition}
Let $W$ be a weak Hopf algebra over $K$. A finite field extension $E/K$ is
called $W$-Galois if there exists a weak Hopf module algebra action
$\alpha\colon W\o_K E\to E$ such that the map
\[
\Phi\colon E\o_L W\to \End_KE\,,\qquad x\o w\mapsto \{y\mapsto x\alpha(w\o y)\}
\]
where $L$ stands for $W\L$, is an isomorphism.
\end{definition}

Let $E/K$ be a finite extension which is $W$-Galois. Below we shall write
$w\triangleright x$ for $\alpha(w\o x)$, $w\in W$, $x\in E$.   

\begin{enumerate}

\item The map $\Phi$ in the above Definition is an $K$-algebra isomorphism from
the smash product $E\rcross W$ to $A$. As a matter of fact, the underlying
$K$-space of the smash product is the tensor product $E\o_LW$ of $L$ modules
where $E_L$ is defined by $x\cdot l:=x(l\triangleright 1)$. So the definition
of $W$-Galois extension just claims that the smash product is isomorphic to
$A$ as an $K$-space. This map is an algebra map as it is obvious from the
multiplication rule of the smash product. Let $\varphi\colon W\to A$ denote
the restriction of this map. 

\item The restriction of the $K$-algebra monomorphism $\varphi\colon W\to A$ to
$L=W\L$ is $\varphi(l)\colon x\mapsto (l\triangleright 1)x$. Therefore
$\varphi(L)\subset E$ and therefore $\varphi$ identifies $L$ with an
intermediate field $K\subset L\subset E$.

\item Let $r\in W\R$. Then $r\triangleright x=x(r\triangleright
1)=(S_W(r)\triangleright 1) x$. Therefore $\varphi(r)=\varphi(S_W(r))\in L$
which, using injectivity of $\varphi$, implies that $W\L=W\R=L$ and $S_W$ acts
as the identity on $L$.

\item It follows that $\cop_W(1_W)$ is a separating idempotent for the
separable algebra $L$ over $K$. By commutativity of $E$ it contains
$\cop_A(1_A)$ as a subprojection, i.e.,
\[
(\varphi\o\varphi)(\cop_W(1_W))\cop_A(1_A)=
\cop_A(1_A)(\varphi\o\varphi)(\cop_W(1_W))=\cop_A(1_A)
\]

\item $\varphi$ is a weak (two sided) morphism of weak bialgebras, i.e.,
\begin{equation} \label{eq: copW vs copA}
\cop_A(1_A)(\varphi\o\varphi)(\cop_W(w))=\cop_A(\varphi(w))\,,\qquad w\in W\,.
\end{equation}
This can be seen as follows.
Upon identifying $A\o_K A$ with $\End_K(E\o_K E)$ the module algebra property
of $_WE$ boils down to
\[
\tau\circ(\varphi\o\varphi)(\cop_W(w))=\tau\circ\cop_A(\varphi(w))\,.
\]
Composing both hand sides with the section $\sigma$ of $\tau$ which is given
by $\sigma(1)=\cop_A(1_A)$ we get precisely the required statement.

\item $W$ is cocommutative. As a matter of fact, module algebra property of
$_WE$ and commutativity of $E$ immediately imply that $(w\1\lact x)(w\2\lact
y)=(w\1\lact y)(w\2\lact x)$ for $x,y\in E$ and $w\in W$. Therefore
$\Phi((w\1\lact y)\rcross w\2)=\Phi((w\2\lact y)\rcross w\1)$ and $\Phi$ being
mono we have equality of the arguments in the smash product. Now the arguments
are images under $v_y\o\id_W$ of $w\1\o_L w\2$ and $w\2\o_L w\1$,
respectively, where $v_y\colon W\to E$, $w\mapsto w\lact y$, is a right
$L$-module map. Since the common kernel of maps $v_y\o\id_W\colon W\o_L W\to
E\o_L W$, $y\in E$, is the kernel of $\varphi\o\id$, we conclude that
\[
w\1\o_L w\2\ =\ w\2\o_L w\1\,,\qquad w\in W\,.
\]
Now use separability of $L/K$ and the isomorphism $W\o_L W\cong (W\o_K
W)\cop_W(1_W)$ to conclude that $w\1\o_K w\2\ =\ w\2\o_K w\1$, $w \in W$.

\item If $n\in\ _WE$ is an invariant then $\varphi(n)$ commutes with
$\varphi(W)$ and since $\End E_K$ is generated by $\varphi(W)$ and the
commutative $E$, it belongs to \break$\Center\End E_K=K$. This proves that
$E^W=K$.

\item $\pi^L_A(\varphi(w))=w\lact 1=\pi^L_W(w)\lact 1=\pi^L_W(w)$, the last
equation identifying $W^L$ with $L\subset E$. Therefore $\eps_A(\varphi(w))=
\eps_A(\pi^L_A(\varphi(w)))=\eps_A(\pi^L_W(w))=\tau(\pi^L_W(w))$. Since
$\eps_W|_L$ is nondegenerate there exist $u\in L_\x$ such that
$\tau(l)=\eps_W(ul)$ for $l\in L$. It follows that
\[
\eps_A(\varphi(w))\ =\ \eps_W(u\pi^L_W(w))=\eps_W(\pi^L_W(uw))
=\eps_W(uw)
\]
thus
\begin{equation} \label{eq: epsW vs epsA}
\eps_W(w)\ =\ \eps_A(\varphi(u^{-1}w))\,,\qquad w\in W\,.
\end{equation}
\end{enumerate}

In particular if $E/K$ is $H$-Galois for some finite dimensional Hopf algebra
$H$ over $K$ then $H$ is embedded into $A$ by a unique weak morphism of
weak Hopf algebras, which is just the restriction of the Galois map. Moreover
$A$ is the crossed product of $E$ with $H$. 

\begin{example}[Greither-Pareigis \cite{Gr-Pa,Mo}]
Let $K=\mathbb{Q}$ the rational field and $E=\mathbb{Q}[x]/(x^4-2)$. Then
$E/K$ is separable but not normal. However it is $H$-Galois for two different
Hopf algebras. One of these Hopf algebras is the commutative and
co-commutative Hopf algebra $H=\mathbb{Q}[c,s]/(c^2+s^2-1,cs)$ with
comultiplication $\cop(c)=c\o c-s\o s$, $\cop(s)=c\o s+s\o c$, counit
$\eps(c)=1$, $\eps(s)=0$ and antipode $S(c)=c$, $S(s)=-s$. The weak left
embedding of $H$ into $A$ is  clear from the presentation of $A$ as 
\begin{equation}
A=\mathbb{Q}\alg\bra c,s,x|c^2+s^2-1,cs,sc,cx-xs,sx+xc,x^4-2\ket
\end{equation}
The general WHA structure of (\ref{WHA A}) can be cast into the form 
\begin{eqnarray}
\cop(1)&=&\frac{1}{4}1\o 1+\frac{1}{8}\sum_{k=1}^3 x^k\o x^{4-k}\\
\cop(x)&=&\frac{1}{4}(x\o 1+1\o x)+\frac{1}{8}(x^3\o x^2+x^2\o x^3)\\
\cop(c)&=&\cop(1)(c\o c-s\o s)\\
\cop(s)&=&\cop(1)(c\o s+s\o c)\\
\eps(c)&=&4\quad\eps(s)=0\quad\eps(x)=0\\
S(c)&=&c\quad S(s)=-s\quad S(x)=x
\end{eqnarray}
\end{example}

\subsection{Galois connection}

Let $E/K$ be separable and let $A$ be its universal weak Hopf algebra. 
Define $\Sub_{Alg/K}(E)$ to be the set of subobjects of $E$ in the category of
$K$-algebras. Also, let $\Sub_{WHA/K}(A)$ be the sub-WHA's of $A$. The latter
means any $K$-subalgebra of $A$ which is closed under comultiplication. 
That is to say we restrict ourselves to strict embeddings of WHA's in the
sense of Definition \ref{def: str/weak mor}.

We can define two order reversing functions (contravariant functors between
preorders)
\[
\Sub_{WHA/K}(A)\rarr{\Fix}\Sub_{Alg/K}(E)\rarr{\Gal}\Sub_{WHA/K}(A)
\]
as follows.
\begin{eqnarray}
\Fix(W)&:=&\{\,x\in E\,|\,a(x)=\pi^L(a)x,\ \forall a\in W\,\}\\
\Gal(F)&:=&\{\,a\in A\,|\,a(xy)=a(x)y,\ \forall x\in E,\,y\in F\,\}
\end{eqnarray}
Then we have the adjointness relations
\begin{equation}
W\subset\Gal(F)\Leftrightarrow F\subset \Fix(W) 
\end{equation}
Therefore the pair $\bra \Fix,\Gal\ket$ is a Galois connection between
sub-WHA's of $A=\End E_K$ and intermediate fields $E\supset F\supset K$.
Moreover this is a half Galois correspondence since every intermediate field
occurs as a fixed field of a sub-WHA,
\begin{equation}
F=\Fix(\Gal(F))\,,\qquad\forall F\in\Sub_{Alg/K}(E)\,.
\end{equation}
A full Galois correspondence would require further analysis of weak
Hopf subalgebras and coideal subalgebras like in \cite{NV2}.

\end{document}